# FREE PRODUCTS OF ABSOLUTE GALOIS GROUPS*

by

Dan Haran

*School of Mathematics, Tel Aviv University*

*Ramat Aviv, Tel Aviv 69978, Israel*

*e-mail: haran@math.tau.ac.il*

and

Moshe Jarden

*School of Mathematics, Tel Aviv University*

*Ramat Aviv, Tel Aviv 69978, Israel*

*e-mail: jarden@math.tau.ac.il*

and

Jochen Koenigsmann

*Fakultät für Mathematik, Universität Konstanz*

*Postfach 5560, 78434 Konstanz, Germany*

*e-mail: koenigsm@informatik.uni-konstanz.de*

Directory: \Jarden\Diary\FreeProd

20 April, 2000

---

* Research supported by the Minkowski Center for Geometry at Tel Aviv University, established by the Minerva Foundation

**Introduction**

Extending [Hei, Thm. 3.2], Efrat and the first author [EfH, Lemma 1.3] prove that if $n$ pro-p groups are absolute Galois groups of fields of the same characteristic, then their free pro-$p$ product is also an absolute Galois group of the same characteristic. Problem 18 of [Jar] asks whether the free profinite product of $n$ absolute Galois groups of fields is an absolute Galois group of a field.

The goal of this note is to give an affirmative answer to the problem. Indeed, we prove a more precise result.

THEOREM:

(a) Suppose $G_1, \ldots, G_n$ are absolute Galois groups of fields. Then their free product $\mathbin{\text{\Large$*$}}_{i=1}^n G_i$ is an absolute Galois group of a field of characteristic 0.

(b) If $G_1, \ldots, G_n$ are absolute Galois groups of fields of a common characteristic $p$, then $\mathbin{\text{\Large$*$}}_{i=1}^n G_i$ is an absolute Galois group of a field of characteristic $p$.

Mel'nikov [Mel, Thm. 1.4], proves the Theorem when $\mathrm{rank}(G_i) \leq \aleph_0$, $i = 1, \ldots, n$. His proof uses a theorem of Geyer [Gey]: Suppose $M$ and $L$ are Henselian fields with respect to rank 1 valuations and both are separable algebraic extensions of a countable field $K$. Then $G(L^\sigma \cap M^\tau) \cong G(L) * G(M)$ for almost all $(\sigma, \tau) \in G(K) \times G(K)$. Here $G(K)$ is the absolute Galois group of $K$ and "almost all" is used in the sense of the Haar measure of $G(K)$. Mel'nikov's proof does not extend to the case of uncountable rank.

Ershov [Ers, Thm. 3] proves part (a) of the theorem by a different method. Our proof simplifies that of Ershov. We take this opportunity to supply proofs to well known results which are not well documented in the literature.

The Theorem is also a consequence of [Pop, Thm. 3.4]. As in the proof of Theorem 3.4, one may use Proposition 2.5 to obtain, in the terminology of [Pop, Ch. 3, §2], a Galois approximation of $\mathbin{\text{\Large$*$}}_{i=1}^n G_i$. Then one map apply [Pop, Thm. 3.4].

# 1. Free products of pro-finite groups

Consider a profinite group $G$ and closed subgroups $G_1, \ldots, G_n$ which generate $G$. We say $G$ is the **free product** of $G_1, \ldots, G_n$ and write $G = \coprod_{i=1}^n G_i$ if the following condition holds:

(1) Given homomorphisms $\eta_i$ of $G_i$ into a profinite group $H$, there is a unique homomorphism $\eta \colon G \to H$ whose restriction to $G_i$ is $\eta_i$, $i = 1, \ldots, n$.

Here is an equivalent condition.

(2) Suppose $\eta_i \colon G_i \to H$, $i = 1, \ldots, n$ are homomorphisms with $H$ finite. Then there is a homomorphism $\eta \colon G \to H$ whose restriction to $G_i$ is $\eta_i$, $i = 1, \ldots, n$.

The construction of $\alpha$ in (2) can be achieved by solving special embedding problems.

Define an **embedding problem** for $(G, G_1, \ldots, G_n)$ to be a structure

(3) $$(\varphi \colon G \to A,\ \psi \colon B \to A,\ B_1, \ldots, B_n)$$

where $\varphi$ and $\psi$ are epimorphisms of profinite groups, $B_1, \ldots, B_n$ are closed subgroups of $B$ which generate $B$, and $\psi$ maps $B_i$ isomorphically onto $\varphi(G_i)$, $i = 1, \ldots, n$. The embedding problem is **finite** if $B$ is finite. A **solution** of (3) is an epimorphism $\gamma \colon G \to B$ with $\psi \circ \gamma = \varphi$ and $\gamma(G_i) = B_i$, $i = 1, \ldots, n$.

LEMMA 1.1: *Let $G$ be a profinite group and $G_1, \ldots, G_n$ closed subgroups which generate $G$. Suppose each finite embedding problem (3) for $(G, G_1, \ldots, G_n)$ has a solution. Then $G = \coprod_{i=1}^n G_i$.*

*Proof:* First, we strengthen the hypothesis of the lemma.

CLAIM A: *If embedding problem (3) has a solution $\gamma$, it is unique.* Indeed, $\psi$ maps $B_i$ bijectively onto $A_i = \varphi(G_i)$. Hence, $\gamma|_{G_i} = (\psi|_{B_i})^{-1} \circ (\varphi|_{G_i})$. This uniquely determines $\gamma|_{G_i}$, $i = 1, \ldots, n$. As $G = \langle G_1, \ldots, G_n \rangle$, $\gamma$ is unique.

CLAIM B: *Each embedding problem (3) with $A$ finite has a solution.* Indeed, there is an inverse system of epimorphisms $B \xrightarrow{\pi^{(j)}} B^{(j)} \xrightarrow{\psi^{(j)}} A$ with finite groups $B^{(j)}$ and $\psi = \psi^{(j)} \circ \pi^{(j)}$ whose inverse limit is $B \xrightarrow{\psi} A$. For each $j$ there is an epimorphism $\gamma^{(j)}\colon G \to B^{(j)}$ with $\psi^{(j)} \circ \gamma^{(j)} = \varphi$ and $\gamma^{(j)}(G_i) = \pi^{(j)}(B_i)$, $i = 1, \ldots, n$. By Claim A, $\gamma^{(j)}$ is unique. Hence, the $\gamma_j$'s are compatible. So, they define an epimorphism $\gamma\colon G \to B$ with $\psi \circ \gamma = \varphi$ and $\gamma(G_i) = B_i$, $i = 1, \ldots, n$.

END OF PROOF: Let $H$ be a finite group and $\eta_i\colon G \to H$ homomorphisms, $i = 1, \ldots, n$. Put $H_i = \eta_i(G_i)$, $i = 1, \ldots, n$. By (2), it suffices to construct a homomorphism $\eta\colon G \to H$ whose restriction to $G_i$ is $\eta_i$, $i = 1, \ldots, n$.

Since $K_i = \mathrm{Ker}(\eta_i)$ is open in $G_i$, there is an open normal subgroup $N$ of $G$, independent of $i$, with $L_i = G_i \cap N \subseteq K_i$, $i = 1, \ldots, n$. Consider the finite group $G/N = \langle G_1 N/N, \ldots, G_n N/N \rangle$ and the free product $B = \ast_{i=1}^{n} G_i/L_i$. Let $\varphi\colon G \to G/N$ be the canonical map $\varphi(g) = gN$. For each $i$ let $\psi_i\colon G_i/L_i \to G_i N/N$ be the canonical isomorphism $\psi_i(gL_i) = gN$. By (2), there is an epimorphism $\psi\colon B \to G/N$ whose restriction to $G_i/L_i$ is $\psi_i$.

Thus, $(\varphi\colon G \to G/N,\ \psi\colon B \to G/N,\ G_1/L_1, \ldots, G_n/L_n)$ is an embedding problem for $(G, G_1, \ldots, G_n)$ with $G/N$ finite. By Claim B, there is an epimorphism $\gamma\colon G \to B$ with $\psi \circ \gamma = \varphi$ and $\gamma(G_i) = G_i/L_i$, $i = 1, \ldots, n$.

For each $i$ let $\bar{\eta}_i\colon G_i/L_i \to H_i$ be the homomorphism with $\bar{\eta}_i(gL_i) = \eta_i(g)$. By (2), there is a homomorphism $\bar{\eta}\colon B \to H$ whose restriction to $G_i/L_i$ is $\bar{\eta}_i$, $i = 1, \ldots, n$.

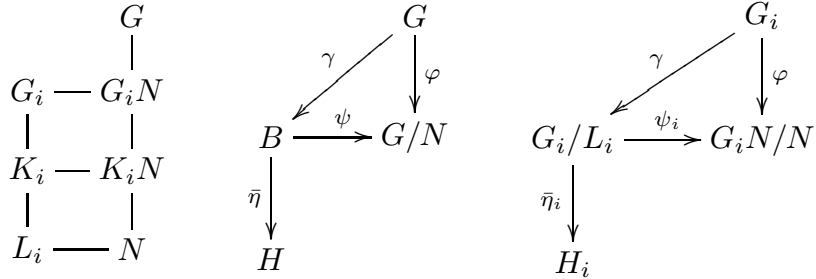

Let $\eta = \gamma \circ \bar{\eta}$. For each $i$ and each $g \in G_i$ we have $\eta(g) = \bar{\eta}_i(\gamma(g)) = \bar{\eta}_i(\psi_i^{-1}(gN)) = \bar{\eta}_i(gL_i) = \eta_i(g)$, as desired. ∎

## 2. Valued fields

The results of this section are well known, although there is some novelty in the presentation. We begin with a brief review of inertia and ramification groups.

Denote the residue field of a valued field $(F,v)$ by $\bar{F}$. For each $x \in F$ with $v(x) \geq 0$ let $\bar{x}$ be the residue of $x$ in $\bar{F}$. Finally, let $F_s$ be the separable algebraic closure of $F$.

Consider a Galois extension $(N,v)/(F,v)$ of Henselian fields. Then $\bar{N}/\bar{F}$ is a normal extension. For each $\sigma \in \mathcal{G}(N/F)$ define $\bar{\sigma} \in \mathrm{Aut}(\bar{N}/\bar{F})$ by the this rule: $\bar{\sigma}\bar{x} = \overline{\sigma x}$ for $x \in N$ with $v(x) \geq 0$. The map $\sigma \mapsto \bar{\sigma}$ is an epimorphism $\rho \colon \mathcal{G}(N/F) \to \mathrm{Aut}(\bar{N}/\bar{F})$ [End, Thm. 19.6]. Its kernel is the **inertia group**:

$$\mathcal{G}_0(N/F) = \{\sigma \in \mathcal{G}(N/F) \mid v(\sigma x - x) > 0 \text{ for each } x \in N \text{ with } v(x) \geq 0\}.$$

Denote the fixed field in $N$ of $\mathcal{G}_0(N/F)$ by $N_0$. Then $\bar{N}_0$ is the maximal separable extension of $\bar{F}$ in $\bar{N}$ [End, Thm. 19.12]. So, $\bar{N}_0/\bar{F}$ is Galois and there is a short exact sequence

$$(1) \qquad 1 \longrightarrow \mathcal{G}_0(N/N_0) \longrightarrow \mathcal{G}(N/F) \xrightarrow{\rho} \mathcal{G}(\bar{N}_0/\bar{F}) \longrightarrow 1.$$

Here we have identified each $\bar{\sigma} \in \mathrm{Aut}(\bar{N}/\bar{F})$ with its restriction to $\bar{N}_0$. In addition, $v(N_0^\times) = v(F^\times)$ [End, Cor. 19.14]. So, $N_0/F$ is an unramified extension.

The **ramification group** of $\mathcal{G}(N/F)$ is

$$\mathcal{G}_1(N/F) = \{\sigma \in \mathcal{G}(N/F) \mid v\left(\frac{\sigma x}{x} - 1\right) > 0 \text{ for each } x \in N\}.$$

It is a normal subgroup of $\mathcal{G}(N/F)$ which is contained in $\mathcal{G}_0(N/F)$ [End, (20.8)]. Denote the fixed field of $\mathcal{G}_1(N/F)$ in $N$ by $N_1$. When $p = \mathrm{char}(\bar{F}) > 0$, $\mathcal{G}_1(N/F)$ is the unique $p$-Sylow subgroup of $\mathcal{G}(N/N_0)$ [End, Thm. 20.18]. When $\mathrm{char}(\bar{F}) = 0$, $\mathcal{G}_1(N/F)$ is trivial. So, in both cases, $\mathrm{char}(\bar{F})$ does not divide $[N_1 : N_0]$.

Suppose now $N = F_s$. Then $N_0 = F_u$ is the **inertia field** and $N_1 = F_r$ is the **ramification field** of $F$. In this case (1) becomes the short exact sequence

$$(2) \qquad 1 \longrightarrow G(F_u) \longrightarrow G(F) \xrightarrow{\rho} G(\bar{F}) \longrightarrow 1.$$

Also, $F \subseteq F_u \subseteq F_r \subseteq F_s$, $F_u/F$ and $F_r/F$ are Galois extensions, $\text{char}(\bar{F}) \nmid [F_r : F_u]$, and $G(F_r)$ is a pro-$p$ group.

Consider now a finite extension $(L, v)/(K, v)$ of Henselian fields. Let $e = e(L/K) = (v(L^\times) : v(K^\times))$ be the **ramification index**. There is a positive integer $d$ such that $[L : K] = de[\bar{L} : \bar{K}]$. If $\text{char}(\bar{K}) = p > 0$, then $d$ is a power of $p$ [Art, p. 62, Thm. 10]. If $\text{char}(\bar{K}) = 0$, then $d = 1$. When $d = 1$ we say $L/K$ is **defectless**. An arbitrary algebraic extension $(M, v)/(K, v)$ is **defectless** if each finite subextension is defectless. This is the case when $\text{char}(\bar{K}) \nmid [L : K]$.

Let $N/F$ be a Galois subextension of $F_u/F$. By (1) and (2), $\bar{N}/\bar{F}$ is Galois and $\rho \colon \mathcal{G}(N/F) \to \mathcal{G}(\bar{N}/\bar{F})$ is an isomorphism. So, $F_u/F$ is defectless. Combine this with the conclusion of the preceding section to conclude: $F_r/F$ is defectless.

LEMMA 2.1: *Let $(F, v)$ be a Henselian valued field. Use the above notation.*
(a) *There is a field $F'$ with $F_u F' = F_r$ and $F_u \cap F' = F$.*
(b) *The short exact sequence $1 \to \mathcal{G}(F_r/F_u) \to \mathcal{G}(F_r/F) \to \mathcal{G}(F_u/F) \to 1$ splits.*

*Proof:* Statement (b) is a Galois theoretic interpretation of (a). So, we prove (a).

Zorn's lemma gives a maximal extension $F'$ of $F$ in $F_r$ with residue field $\bar{F}$. For each prime number $l \neq \text{char}(\bar{F})$ the value group of $F'$ is $l$-divisible Otherwise, there is $a \in F'$ with $v(a) \notin lv((F')^\times)$. Hence, $F'(\sqrt[l]{a})/F'$ is an extension of degree $l$ and ramification index $l$. Recall: $G(F_r)$ is a pro-$p$ group if $\text{char}(\bar{F}) = p > 0$ and trivial if $\text{char}(\bar{F}) = 0$. So, $F'(\sqrt[l]{a}) \subseteq F_r$. The residue field of $F'(\sqrt[l]{a})$ coincides with $\bar{F'} = \bar{F}$, in contrast to the maximality of $F'$.

By the discussion preceding Lemma 2.1, $F_u \cap F' = F$. Let $E = F_u F'$. Consider a prime number $l \neq \text{char}(\bar{F})$. Then $v(E^\times)$ is contained in the $l$-divisible hull of $v(F')$. As $v((F')^\times)$ is $l$-divisible, so is $v(E^\times)$. In addition $\bar{E} = \bar{F}_u = \bar{F}_s$. As $[F_r : E]$ divides $[F_r : F_u]$, it is not divisible by $\text{char}(K)$. So, $F_r/E$ is defectless. Conclude that $E = F_r$.
∎

LEMMA 2.2 (Kuhlmann-Pank-Roquette [KPR, Thm. 2.2]): *Let $(F, v)$ be a Henselian field.*
(a) *There is a field $F'$ with $F_r \cap F' = F$ and $F_r F' = F_s$.*

(b) *The short sequence* $1 \to G(F_r) \to G(F) \to \mathcal{G}(F_r/F) \to 1$ *splits.*

Proof: Let $p = \text{char}(\bar{F})$. If $p = 0$, then $F_r = F_s$ and we may take $F' = F$. Suppose $p \neq 0$.

By (2), $\mathcal{G}(F_u/F) \cong G(\bar{F})$. By Witt, the $p$-Sylow subgroups of $G(\bar{F})$ are free [Rib, p. 256, Thm. 3.3]. Hence, so are the $p$-Sylow subgroups of $\mathcal{G}(F_u/F)$. As $p \nmid [F_r : F_u]$, restriction $\mathcal{G}(F_r/F) \to \mathcal{G}(F_u/F)$ maps each $p$-Sylow subgroup of $\mathcal{G}(F_r/F)$ isomorphically onto a $p$-Sylow subgroup of $\mathcal{G}(F_u/F)$. So, each $p$-Sylow subgroup of $\mathcal{G}(F_r/F)$ is free. Thus $\text{cd}_p(\mathcal{G}(F_r/F)) = 1$ [Rib, p. 207, Cor. 2.2]. As $G(F_r)$ is a pro-$p$ group, the short sequence in (b) splits [Rib, p. 211, Prop. 3.1(iii)']. ∎

PROPOSITION 2.3: *Let $(F, v)$ be a valued field.*

(a) *Suppose $(F, v)$ is Henselian. Then the epimorphism $\rho\colon G(F) \to G(\bar{F})$ induced by reduction at $v$ has a section.*

(b) *Each subgroup of $G(\bar{F})$ is isomorphic to a subgroup of $G(F)$.*

Proof of (a): The map $\rho$ decomposes as $G(F) \xrightarrow{\text{res}} \mathcal{G}(F_r/F) \xrightarrow{\text{res}} \mathcal{G}(F_u/F) \xrightarrow{\bar{\rho}} G(\bar{F})$. The map $\bar{\rho}$ which is also induced by reduction is an isomorphism (by (1)). By Lemmas 2.1 and 2.2, each of the restriction maps splits. Hence $\rho$ splits.

Proof of (b): Let $(F', v)$ be the Henselization of $(F, v)$. Then $\overline{F'} = \bar{F}$. By (a), each subgroup of $G(\bar{F})$ is isomorphic to a subgroup of $G(F')$, hence of $G(F)$. ∎

PROPOSITION 2.4: *Let $F/K$ be an extension of fields and $v$ a valuation of $F$. Suppose $v$ is trivial on $K$ and $\bar{F} = K$. Then* res$\colon G(F) \to G(K)$ *is an epimorphism which has a section. This conclusion holds in particular when $F/K$ is a purely transcendental extension.*

Proof: Replace $(F, v)$ by its Henselian closure, if necessary, to assume $(F, v)$ is Henselian. For each $a \in K$ and each $\sigma \in G(K)$ we have, $\sigma a = \overline{\sigma a} = \bar{\sigma}\bar{a} = \bar{\sigma}a = \rho(\sigma)a$. So, restriction coincides with the map $\rho$ induced by reduction. Now apply Proposition 2.3(a).

When $F/K$ is a purely transcendental extension, $F$ has a valuation with residue field $K$. This is evident when $F = K(t)$ and $t$ is transcendental. The general case follows from the special case by transfinite induction and using composition of valuations. So, res: $G(F) \to G(K)$ has a section. ∎

The following Proposition gives more details to a result of Efrat [Efr. Prop. 4.7].

PROPOSITION 2.5: *Let $K$ be a field, $E_0$ its prime field, and $T$ a set of variables with $\mathrm{card}(T) \geq \mathrm{trans.deg}(K/E_0)$. Let $F_0$ be either $E_0$ or $\mathbb{Q}$. Then there is a field $L$, algebraic over $F_0(T)$, with $G(L) \cong G(K)$.*

*Proof:* There is a unique place $\varphi_0 \colon F_0 \to E_0 \cup \{\infty\}$. Choose a transcendence base $\bar{T}$ for $K/E_0$. By assumption, $\mathrm{card}(\bar{T}) \leq \mathrm{card}(T)$. Choose a surjective map $\varphi_1 \colon T \to \bar{T}$. Let $F_1 = F_0(T)$ and $E_1 = E_0(\bar{T})$. Extend $\varphi_0$ and $\varphi_1$ to a place $\varphi \colon F_1 \to E_1 \cup \{\infty\}$. Denote the corresponding valuation by $v$. Corollary 2.3(b) gives the desired field $L$. ∎

## 3. Embedding problems

Starting from a field $K$ and separable algebraic extensions $K_1, \ldots, K_n$ with $\bigcap_{i=1}^{n} K_i = K$, we construct extensions $L, L_1, \ldots, L_n$ such that $\bigcap_{i=1}^{n} L_i = L$, res: $G(L_i) \to G(L_i)$ is bijective, $i = 1, \ldots, n$, and each finite embedding problem for $(G(L), G(L_1), \ldots, G(L_n))$ is solvable.

The special case of the following lemma where $H = G$ is [FrJ, Lemma 24.44]. The general case essentially appears in [HaJ, Part C of the proof of Proposition 14.1].

LEMMA 3.1: *Let $L/K$ be a finite Galois extension. Let $\psi\colon B \to \mathcal{G}(L/K)$ be an epimorphism of finite groups. Then there exists a finitely generated regular extension $E$ of $K$ and a finite Galois extension $F$ of $E$ containing $L$ such that $B = \mathcal{G}(F/E)$ and $\psi$ is the restriction $\mathrm{res}_{F/L}\colon \mathcal{G}(F/E) \to \mathcal{G}(L/K)$.*

*Moreover, let $K \subseteq L_0 \subseteq L$ be a field and $B_0$ a subgroup of $B$ which $\psi$ maps isomorphically onto $\mathcal{G}(L/L_0)$. Then the fixed field $F_0$ of $B_0$ in $F$ is a purely transcendental extension of $L_0$.*

*Proof:* Let $x^\beta$, $\beta \in B$, be algebraically independent elements over $K$. Define an action of $B$ on $F = L(x^\beta \mid \beta \in B)$ by $(x^\beta)^{\beta'} = x^{\beta\beta'}$ and $a^{\beta'} = a^{\psi(\beta')}$ for $a \in L$. Denote the fixed field of $B$ in $F$ by $E$. Then $F/K$ is a finitely generated separable extension. By [Lan, p. 64, Prop. 6], $E/K$ is also a finitely generated separable extension. Also, res: $\mathcal{G}(F/E) \to \mathcal{G}(L/K)$ coincides with $\psi\colon B \to \mathcal{G}(L/K)$. So, $E \cap \tilde{K} = E \cap F \cap \tilde{K} = E \cap L = K$. Conclude: $E/K$ is regular.

Choose a set of representatives $R$ for the left cosets of $B$ modulo $B_0$. Let $w_1, \ldots, w_m$ be a basis for $L/L_0$. By assumption, $m = |B_0|$. Consider $\rho \in R$. Put

$$(1) \qquad t_{\rho j} = \sum_{\beta \in B_0} w_j^\beta x^{\rho\beta}, \qquad j = 1, \ldots, m.$$

Since $\det(w_j^\beta) \neq 0$, each $x^{\rho\beta}$ is a linear combination of $t_{\rho j}$ with coefficients in $L$. Put $\mathbf{t} = (t_{\rho j} \mid \rho \in R, \ j = 1, \ldots, m)$, $\mathbf{x} = (x^\beta \mid \beta \in B)$, and $n = |B|$. Both tuples contain exactly $n$ elements and $L(\mathbf{t}) = L(\mathbf{x}) = F$. So, $L_0(\mathbf{t})$ is a purely transcendental extension of $L_0$.

Each $t_{\rho j}$ is fixed by $H$. So, $L_0(\mathbf{t}) \subseteq F_0$. Moreover, $m = [L : L_0] = [L(\mathbf{t}) : L_0(\mathbf{t})] \geq [F : F_0] = |H| = m$. Conclude: $F_0 = L_0(\mathbf{t})$ and $F_0/L_0$ is purely transcendental. ∎

LEMMA 3.2: Let $K, K_1, \ldots, K_n$ be fields with $K_i/K$ separable algebraic, $i = 1, \ldots, n$, and $\bigcap_{i=1}^n K_i = K$. Let $(\varphi\colon G(K) \to A, \psi\colon B \to A, B_1, \ldots, B_n)$ be a finite embedding problem for $(G(K), G(K_1), \ldots, G(K_n))$. Then there are fields $L, L_1, \ldots, L_n$ with these properties:

(2a) $L/K$ has a positive finite transcendence degree.

(2b) $L_i/L$ is separable algebraic, $i = 1, \ldots, n$, and $\bigcap_{i=1}^n L_i = L$.

(2c) $K_i \subseteq L_i$ and res: $G(L_i) \to G(K_i)$ is a bijection, $i = 1, \ldots, n$.

(2d) There is a homomorphism $\gamma\colon G(L) \to B$ with $\psi \circ \gamma = \varphi \circ \mathrm{res}_{L_s/K_s}$ and $\gamma(G(L_i)) = B_i$, $i = 1, \ldots, n$.

Proof: Let $N$ be the fixed field of $\mathrm{Ker}(\varphi)$ in $K_s$. It is a finite Galois extension of $K$ and $\mathcal{G}(N/K) \cong A$. Assume without loss $A = \mathcal{G}(N/K)$ and $\varphi = \mathrm{res}_{K_s/N}$. Let $N_i = K_i \cap N$, and $A_i = \mathcal{G}(N/N_i)$. Then $A_i = \varphi(G(K_i))$, $i = 1, \ldots, n$.

By Lemma 3.1, there is a finitely generated regular extension $E$ of $K$ and a finite Galois extension $F$ of $E$ containing $N$ such that $B = \mathcal{G}(F/E)$ and $\psi = \mathrm{res}_{F/N}$. Moreover, for each $1 \leq i \leq n$, the fixed field $F_i$ of $B_i$ in $F$ is a purely transcendental extension of $N_i$. Corollary 2.4 gives a section $\lambda_i\colon G(N_i) \to G(F_i)$ to the map res: $G(F_i) \to G(N_i)$. Denote the fixed field of $\lambda_i(G(K_i))$ in $E_s$ by $L_i$. Then res: $G(L_i) \to G(K_i)$ is a bijection.

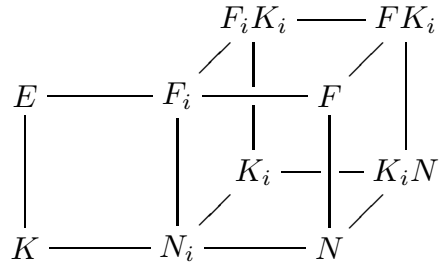

Finally, let $L = \bigcap_{i=1}^n L_i$ and $\gamma\colon G(L) \to B$ be restriction. Then $\gamma(G(L))$ contains $B_i = \gamma(G(L_i))$, $i = 1, \ldots, n$. Since $B = \langle B_1, \ldots, B_n \rangle$, $\gamma$ is surjective. As all maps are now restrictions, $\psi \circ \gamma = \varphi \circ \mathrm{res}_{L_s/K_s}$, as desired. ∎

THEOREM 3.3: *Let* $K, K_1, \ldots, K_n$ *be fields with* $K_i/K$ *separable algebraic*, $i = 1, \ldots, n$, *and* $\bigcap_{i=1}^n K_i = K$. *Then there are fields* $F, F_1, \ldots, F_n$ *with these properties*:

(3a) $\operatorname{card}(F) = \max(\aleph_0, \operatorname{card}(K))$.

(3b) $F_i/F$ *is separable algebraic*, $i = 1, \ldots, n$, *and* $\bigcap_{i=1}^n F_i = F$.

(3c) $K_i \subseteq F_i$ *and* $\operatorname{res}\colon G(F_i) \to G(K_i)$ *is a bijection*, $i = 1, \ldots, n$.

(3d) $G(F) = \boxast_{i=1}^n G(F_i)$.

*Proof:* Put $m = \max(\aleph_0, \operatorname{card}(K))$. Let

$$(\varphi_\alpha\colon G(K) \to A_\alpha,\ \psi_\alpha\colon B_\alpha \to A_\alpha,\ B_{\alpha,1}, \ldots, B_{\alpha,n}), \quad \alpha < m,$$

be a well-ordering of all embedding problems for $(G(K), G(K_1), \ldots, G(K_n))$. Multiple application of Lemma 3.2 and transfinite induction gives for each ordinal $\alpha < m$ fields $(L_\alpha, L_{\alpha,1}, \ldots, L_{\alpha,n})$ with $L_0 = K$ and $L_{0,i} = K_i$, $i = 1, \ldots, n$. Moreover, for all $\alpha < \beta < m$ the following holds:

(4a) $L_{\alpha+1}/L_\alpha$ has a positive finite transcendence degree.

(4b) $L_{\alpha,i}/L_\alpha$ is separable algebraic, $i = 1, \ldots, n$, and $\bigcap_{i=1}^n L_{\alpha,i} = L_\alpha$.

(4c) $L_{\alpha,i} \subseteq L_{\beta,i}$ and $\operatorname{res}\colon G(L_{\beta,i}) \to G(L_{\alpha,i})$ is a bijection, $i = 1, \ldots, n$.

(4d) There is an epimorphism $\gamma_{\alpha+1}\colon G(L_{\alpha+1}) \to B_\alpha$ such that $\psi_\alpha \circ \gamma_\alpha = \varphi_{\alpha+1} \circ \operatorname{res}_{L_{\alpha+1,s}/K_s}$ and $\gamma_\alpha(G(L_{\alpha+1,i})) = B_{\alpha,i}$, $i = 1, \ldots, n$.

(4e) $L_\beta = \bigcup_{\alpha<\beta} L_\alpha$ and $L_{\beta,i} = \bigcup_{\alpha<\beta} L_{\alpha,i}$, $i = 1, \ldots, n$, when $\beta$ is a limit ordinal.

Indeed, suppose all objects with index $\alpha$ have been constructed. As $K = \bigcap_{i=1}^n K_i$, $L_\alpha = \bigcap_{i=1}^n L_{\alpha,i}$, and $\operatorname{res}(G(L_{\alpha,i})) = G(K_i)$, the map $\operatorname{res} \circ \varphi\colon G(L_\alpha) \to A_\alpha$ is surjective. We may therefore apply Lemma 3.2 to $L_\alpha, L_{\alpha,1}, \ldots, L_{\alpha,n}, \varphi_\alpha \circ \operatorname{res}_{L_{\alpha,s}/K_s}$ replacing $K, K_1, \ldots, K_n, \varphi_\alpha$ and construct the objects with index $\alpha+1$ that satisfy (4a) and (4d).

Let $M_1 = \bigcup_{\alpha<m} L_\alpha$, $M_{1,i} = \bigcup_{\alpha<m} L_{\alpha,i}$, $i = 1, \ldots, n$. Then the following holds:

(5a) $\operatorname{card}(M_1) = m$.

(5b) $M_{1,i}/M_1$ is separable algebraic, $i = 1, \ldots, n$, and $\bigcap_{i=1}^n M_{1,i} = M_1$.

(5c) $K_i \subseteq M_{1,i}$ and $\operatorname{res}\colon G(M_{1,i}) \to G(K_i)$ is a bijection, $i = 1, \ldots, n$.

(5d) For each $\alpha < m$ there is an epimorphism $\gamma_{1,\alpha}\colon G(M_1) \to B_\alpha$ with $\psi_\alpha \circ \gamma_{1,\alpha} = \varphi_\alpha \circ \operatorname{res}_{M_{1,s}/K_s}$ and $\gamma_{\alpha,1}(G(M_{1,i})) = B_{\alpha,i}$, $i = 1, \ldots, n$.

Finally use ordinary induction to construct an ascending sequence of fields $K = M_0 \subset M_1 \subset M_2 \subset \cdots$ such that $M_j$ satisfies (5) with $M_{j-1}, M_j$ replacing $K, M_1$. Let $F = \bigcup_{j=1}^{\infty} M_j$, $F_i = \bigcup_{j=1}^{\infty} M_{j,i}$, $i = 1, \ldots, n$. Then (3a) and (3b) hold. In particular, $G(F) = \langle G(F_1), \ldots, G(F_n) \rangle$. Moreover, for each $i$ and $j$ the map res: $G(F_i) \to G(M_{j,i})$ is bijective and the map res: $G(F) \to G(M_j)$ is surjective. Each finite embedding problem $\mathcal{E}$ for $(G(F), G(F_1), \ldots, G(F_n))$ factors through an embedding problem for $(G(M_j), G(M_{j,1}), \ldots, G(M_{j,n}))$ for some $j$. The latter gives an embedding problem for $G(M_{j+1}), G(M_{j+1,1}), \ldots, G(M_{j+1,n})$ which is solvable. So, $\mathcal{E}$ is solvable. Conclude from Lemma 1.1 that $G(F) = \mathop{*}\limits_{i=1}^{n} G(F_i)$. ∎

We are now ready to prove the main result.

THEOREM 3.4: *Let $K_1, \ldots, K_n$ be fields. Then there exists a field $F$ of characteristic $0$ with $G(F) \cong \mathop{*}\limits_{i=1}^{n} G(K_i)$. When all $K_i$ have a common characteristic $p$, there exists a field $F$ with* char$(F) = p$ *and* $G(F) \cong \mathop{*}\limits_{i=1}^{n} G(K_i)$.

*Proof:* Choose an infinite set of variables $T$ of cardinality at least as the transcendence degree of $K_i$ over its prime field, $i = 1, \ldots, n$. Put $K = \mathbb{Q}(T)$ or $K = \mathbb{F}_p(T)$, if char$(K_i) = p > 0$ for all $i$. For each $i$, Proposition 2.5 gives an algebraic extension $L_i$ of $K$ with $G(L_i) \cong G(K_i)$. Then Theorem 3.3 gives an extension $F$ of $K$ with $G(F) \cong \mathop{*}\limits_{i=1}^{n} G(K_i)$. ∎

Remark 3.5: *Free products in the category of pro-$p$ groups.* Let $G_1, \ldots, G_n$ be pro-$p$ groups for some prime number $p$. Denote their free product in the category of all profinite groups (resp. pro-$p$ groups) by $G^*$ (resp. $G^{*,p}$). Each of these groups is generated by $G_1, \ldots, G_n$. Denote the epimorphism of $G^* \to G^{*,p}$ whose restriction to each $G_i$ is the identity map by $\alpha$. Choose a $p$-Sylow subgroup $P$ of $G^*$. Then $\alpha(P) = G^{*,p}$.

For each $i$ there is $a_i \in G^*$ with $G_i^{a_i} \leq P$. Choose $b_i \in P$ with $\alpha(b_i) = \alpha(a_i)$. Then $G_i^{a_i b_i^{-1}} \leq P$. The map $g \mapsto g^{a_i b_i^{-1}}$ defines an embedding of $G_i$, considered as a subgroup of $G^{*,p}$, into $P$. Hence, there is a homomorphism $\alpha' \colon G^{*,p} \to P$ such that $\alpha'(g) = g^{a_i b_i^{-1}}$ and therefore $\alpha(\alpha'(g)) = g^{\alpha(a_i)\alpha(b_i)^{-1}} = g$ for each $i$ and each $g \in G_i$. So, $\alpha \circ \alpha'$ is the identity map on $G^{*,p}$. Conclude: $\alpha' \colon G^{*,p} \to G^*$ is an embedding.

Suppose now each $G_i$ is an absolute Galois group. Theorem 3.4 gives a field $K$

with $G(K) \cong G^*$. By the preceding paragraph, $K$ has an extension $L$ with $G(L) \cong G^{*,p}$. Thus, Theorem 3.4 implies Lemma 1.3 of [EfH] mentioned in the introduction. ∎

27 May, 2000